\documentclass[twoside,12pt]{amsart}

\usepackage{amssymb,amsthm,amsfonts,amsmath,pifont}
\usepackage{graphicx}
\usepackage{pstricks,pst-node, pst-text, pst-plot}
\usepackage{color}
\usepackage{enumerate}
\usepackage{wrapfig}
\usepackage[all]{xy}

\title{On the impossibility of constructing a triangle given its internal bisectors}

\author{A. Caminha$^1$}
\address{$^1$Departamento de Matem\'atica, Universidade Federal do Cear\'a, Fortaleza,
Cear\'a, Brazil. 60455-760}
\email{caminha@mat.ufc.br}

\author{A. Maia$^2$}
\address{$^1$Departamento de Matem\'atica, Universidade Federal do Cear\'a, Fortaleza,
Cear\'a, Brazil. 60455-760}
\email{caminha@mat.ufc.br}

\subjclass[2010]{Primary 12F05; Secondary 12-01}

\keywords{Abstract Algebra; Greek problems}

\begin{document}
 \maketitle
 
\begin{abstract}
We give a simple proof to the fact that it is impossible to use straightedge and compass to construct a triangle given the lengths of its internal bisectors, even if the triangle is isosceles.
\end{abstract}

\vspace{1cm}

\begin{figure}[h]\begin{center}
\begin{pspicture}(0,-1)(2,3.5)
\pspolygon(0,0)(2,0)(1,3)
\uput[90](1,3){$A$}\qdisk(1,3){1.2pt}
\uput[-20](1.9,-.05){$C$}\qdisk(2,0){1.2pt}
\uput[200](0.1,-.05){$B$}\qdisk(0,0){1.2pt}

\psline(1,3)(1,0)
\uput[270](1,.05){$M$}\qdisk(1,0){1.2pt}

\psline(0,0)(1.55,1.3)
\uput[80](1.75,1.1){$P$}\qdisk(1.56,1.3){1.2pt}

\psarc(2,0){.25}{110}{180}
\uput[145](1.95,.05){$\theta$}

\psline{|-|}(0,-.5)(2,-.5)
\uput[270](1,-.4){$b$}


\psline{|-|}(.67,3.14)(-.45,0.14)
\uput[150](.11,1.57){$l$}

\end{pspicture}
\caption{\small an isosceles triangle and two of its internal bisectors.}\label{fig:the triangle}
\end{center}\end{figure}
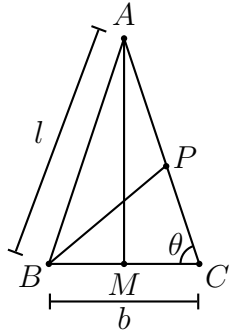

We consider a triangle $ABC$, isosceles of basis $BC$, and let $AM$ and $AP$ be internal bisectors of it (cf. Figure \ref{fig:the triangle}). We assume the lengths $p$ of $AP$ and $q$ of $AM$ to be known.

We let $AB=AC=l$ and $BC=b$. Since $AM$ is also height and median of $ABC$, we have $\cos\theta=\frac{b}{2l}$; also, Pythagoras' theorem applied to triangle $ACM$ gives $l^2-\frac{b^2}{4}=q^2$ or, which is the same,
\begin{equation}\label{eq:auxiliar 1}
(2l+b)(2l-b)=4q^2.
\end{equation}

With respect to the internal bisector AP, the interior angle bisector theorem furnishes $\frac{AP}{CP}=\frac{AB}{BC}$ or, letting $x=CP$, $\frac{l-x}{x}=\frac{l}{b}$. Solving for $x$, we obtain $x=\frac{bl}{b+l}$, and applying the cosine law to triangle $BPC$, we get
\begin{equation}\label{eq:auxiliar 2}
 \begin{split}
p^2&\,=b^2+x^2-2bx\cos\theta\\
&\,=b^2+\Big(\frac{bl}{b+l}\Big)^2-2b\Big(\frac{bl}{b+l}\Big)\cdot\frac{b}{2l}\\
&\,=\frac{b^2l(b+2l)}{(b+l)^2}.
 \end{split}
\end{equation}

From \eqref{eq:auxiliar 1} and \eqref{eq:auxiliar 2}, it comes that
$$\frac{b^2l}{2l-b}=\frac{p^2(b+l)^2}{4q^2}.$$
By cross-multiplying, expanding and performing some elementary algebra, we arrive at the equality $2p^2l^3+3p^2bl^2-4q^2b^2l-p^2b^3=0$
and, upon division by $b^3$,
$$2p^2\Big(\frac{l}{b}\Big)^3+3p^2\Big(\frac{l}{b}\Big)^2-4q^2\Big(\frac{l}{b}\Big)-p^2=0.$$
Assuming, without loss of generality, that $p=1$, we conclude that $\frac{l}{b}$ is a root of the third degree polynomial
$$f(X)=2X^3+3X^2-4q^2X-1.$$

Now, for the sake of contradiction, suppose that one can use straightedge and compass to construct $ABC$, knowing the lengths $p=1$ and $q$. By recalling the usual analysis of the classical Greek construction problems (cf. \cite{Hadlock:00}, for instance), this means that there exists a finite sequence of elementary constructions that allows us to obtain the length $\frac{lp}{b}=\frac{l}{b}$ (since we are assuming $p=1$). Yet in another way, $\frac{l}{b}$ is constructible from $\mathbb Q(p,q)=\mathbb Q(q)$, so that (by arguing again as in the analysis of the Greek problems) the degree $[\mathbb Q(q)(l/b):\mathbb Q(q)]$ must be a power of $2$. However, if we show that $f$ is irreducible in $\mathbb Q(q)[X]$, then
$$[\mathbb Q(q)(l/b):\mathbb Q(q)]=3,$$
which will be a contradiction.

We are left to establishing the irreducibility of $f$ in $\mathbb Q(q)[X]$, at least for some $q>0$. To this end, from now on we take $q>0$ to be transcendental. Then, $\mathbb Q[q]\simeq\mathbb Q[X]$ is a UFD, and Gauss' theorem assures that it suffices to show that $f$ is irreducible in $\mathbb Q[q][X]$. If this is not so, then $f$ has a root $\alpha\in\mathbb Q(q)$, say $\alpha=\frac{g(q)}{h(q)}$, for some $g,h\in\mathbb Q[X]\setminus\{0\}$. 

Applying the searching criterion for roots of $f$ belonging to the field of fractions of the UFD $\mathbb Q[q]$ assures that $g(q)\mid 1$ in $\mathbb Q[q]$. Therefore, we can assume $g(q)=1$, so that $\alpha=\frac{1}{h(q)}$ and $f(\alpha)=0$ gives 
$$2+3h(q)-4q^2h(q)^2=h(q)^3.$$ 
This is the same as 
$$h(X)^2(4X^2+h(X))=3h(X)+2,$$
so that $\partial h=0$, say $h(X)=c$ for some $c\in\mathbb Q$. But this gives
$$c^2(4X^2+c)=3c+2,$$
which is clearly impossible.

\end{document}